\documentclass{article}
\usepackage{arxiv}
\usepackage{graphicx}

\usepackage[utf8]{inputenc} 
\usepackage[T1]{fontenc}    
\usepackage{hyperref}       
\usepackage{url}            
\usepackage{booktabs}       
\usepackage{amsmath,amsthm,amssymb,amsfonts}       
\usepackage{nicefrac}       
\usepackage{microtype}      
\usepackage{bm}
\usepackage{todonotes}
\setuptodonotes{tickmarkheight=0.2cm}

\newcommand{\half}{\frac{1}{2}}
\newcommand{\emh}{{e-\half}}
\newcommand{\eph}{{e+\half}}

\newcommand{\myvec}[1]{\bm{#1}}
\newcommand{\mymat}[1]{\bm{\mathcal{#1}}}

\newcommand{\D}{\mathcal{D}}

\newcommand{\Dm}{\mymat{D}}
\newcommand{\av}{\myvec{a}}
\newcommand{\bv}{\myvec{b}}
\newcommand{\cv}{\myvec{c}}
\newcommand{\rhov}{\myvec{\rho}}
\newcommand{\uv}{\myvec{u}}

\newcommand{\Ev}{\myvec{E}}
\newcommand{\mv}{\myvec{m}}
\newcommand{\pv}{\myvec{p}}

\newcommand{\Mm}{\mymat{M}}
\newcommand{\Bm}{\mymat{B}}
\newcommand{\Rm}{\mymat{R}}

\newcommand{\ro}{r}

\newcommand{\pen}{\mathcal{E}}

\newcommand{\SAT}{\textrm{SAT}}

\newcommand{\avg}[1]{\{#1\}}
\newcommand{\ud}{\textrm{d}}
\newcommand{\od}[2]{\frac{\ud #1}{\ud #2}}
\newcommand{\pd}[2]{\frac{\partial #1}{\partial #2}}
\newcommand{\norm}[1]{\| #1 \|}
\newcommand{\keep}{KEEP-PE }

\title{On a linear stability issue of split form schemes for compressible flows}

\author{
  Vikram Singh
  \\  The Ocean in the Earth System\\
  Max-Planck Institute for Meteorolgy\\
  Hamburg - 20146, Germany \\
  \texttt{vikram.singh@mpimet.mpg.de} \\
   \And
 Praveen Chandrashekar \\
  Center for Applicable Mathematics\\
  Tata Institute of Fundamental Research\\
  Bangalore -- 560065, India\\
  \texttt{praveen@math.tifrbng.res.in} \\
}

\begin{document}

\maketitle

\begin{abstract}
Split form schemes for Euler and Navier-Stokes equations are useful for computation of turbulent flows due to their better robustness. This is because they satisfy additional conservation properties of the  governing equations like kinetic energy preservation leading to a reduction in aliasing errors at high orders.  Recently, linear stability issues have been pointed out  for these schemes for a density wave problem and we investigate this behaviour for some standard split forms. By deriving linearized equations of split form schemes,  we show that most existing schemes do not satisfy a perturbation energy equation that holds at the continuous level.   A simple modification to the energy flux of some existing schemes is shown to yield a scheme that is consistent with the energy perturbation equation. Numerical tests are given using a discontinuous Galerkin method to demonstrate these results.
\end{abstract}
\section{Introduction}

High order methods are essential for large eddy simulation~(LES)  and direct numerical simulation~(DNS) of turbulent flow problems but their stability is a major issue, especially in the under-resolved cases.   High order methods suffer from aliasing instabilities that can lead to solver blow-up~\cite{MOURA2017, Pirozzoli2010, Rojas2021, SINGH2020}, therefore improving robustness of these methods is an active area of research. One approach for improving robustness that has become increasingly popular is  carrying over stability properties of the continuous PDE to the discrete setup, using techniques such as Summation-By-Parts operators (SBP)~\cite{SVARD201417}.  For example, the compressible Euler equations satisfy  kinetic energy preservation (KEP)~\cite{Jameson2008} and entropy preservation~\cite{LeFloch2000} properties under some conditions.  Many different methods have been developed to extend these properties to the discrete setup using so called split forms, which can increase robustness \cite{KENNEDY2008, MORINISHI2010276, Blaisdell1996, Jameson2008, Ducros2000, COPPOLA201986, FISHER2013353}.  However, recent studies \cite{gassner2020stability, ranocha2020preventing} have shown  linear stability issues with this approach. A simple solution of the Euler equations which involves a density wave propagating with constant speed and constant pressure was considered.  The linearization of certain entropy conservative schemes around this density wave solution was shown to exhibit unstable eigenvalues while the standard central scheme does not have this issue\footnote{If the semi-discrete scheme is $\od{U}{t}=R(U)$, its linearization around some state $\bar{U}$ is given by $\od{U'}{t} = R'(\bar U) U'$. For stability, all the eigenvalues of $R'(\bar U)$ must have negative real parts.}. Many split form schemes which are designed to be kinetic energy preserving also suffer from linear instability in the sense of eigenvalues.  


In this study, we analyze several split form schemes for their linear stability for the density wave problem.  First, we derive linearized Euler equations around the density wave solution considered as the base solution, and show that small perturbations around the base solution are not strictly bounded by the initial data but can admit some growth. Since the base solution has small values of density, the perturbations can drive it to still smaller values, possibly negative, and make the solution to be unstable in a non-linear sense, if the underlying schemes are not monotone. While there does not exist an energy norm  for the perturbations which is non-increasing, there exists a reduced energy that is conserved, but it does not bound the density perturbations.

We also linearize the split form schemes around the density wave solution, and investigate their stability in terms of conserving the reduced energy. We find that even the central scheme does not satisfy reduced energy conservation property present in the linearized PDE even though it is stable in the sense of eigenvalues.  We propose a modified KEP scheme that mimics the energy principles of the linearized PDE. Finally, we present numerical results using nodal discontinuous Galerkin methods and compare some split form schemes in terms of their stability and ability to conserve quantities of interest.
\section{DG split form schemes}
\label{sec:dg}
Consider the Euler equations in conservation form which can be written as
\begin{equation}
U_t + F_x + G_y + H_z = 0
\label{eq:3deul}
\end{equation}
where $U$ is the vector of conserved variables and $F,G,H$ are the Cartesian components of the flux vectors given by
\begin{equation}
U = \begin{bmatrix}
\rho \\ \rho u \\ \rho v \\ \rho w \\ E
\end{bmatrix}, \qquad 
F = \begin{bmatrix}
\rho u \\ p + \rho u^2 \\ \rho u v \\ \rho u w \\ (E+p)u 
\end{bmatrix}, \qquad G = \begin{bmatrix}
\rho v \\ \rho u v \\ p + \rho v^2 \\ \rho wv \\ (E+p)v 
\end{bmatrix}, \qquad H = \begin{bmatrix}
\rho w \\ \rho u w \\ \rho v w \\ p + \rho  w^2 \\ (E+p)w 
\end{bmatrix}
\label{eq:UFGH}
\end{equation}
The subscripts $t, x, y, z$ in~\eqref{eq:3deul} refer to the time derivative and spatial derivatives respectively. In~\eqref{eq:UFGH}, $\rho,p$ are the density and pressure, $u,v,w$ are the Cartesian components of the fluid velocity, $E=p/(\gamma-1) + \rho(u^2 + v^2 + w^2)/2$ is the total energy per unit volume, and $\gamma>1$ is the ratio of specific heats at constant pressure and volume.

We consider the  discontinuous Galerkin (DG) schemes  using nodal Lagrange polynomials collocated at the Gauss-Lobatto-Legendre (GLL) nodes to discretize the compressible Euler equations. In the following, we describe the method in 1-D only, and the extension to multiple dimensions can be made using a tensor product approach.  Let $N$ be the degree of polynomials used to approximate the solution; each component of the solution inside the $e$'th element is written as
\[
U^e = \sum_{i=0}^N U^e_i \ell_i(\xi)
\]
where $\ell_i(\xi)$ are the nodal Lagrange polynomials based on $(N+1)$-GLL nodes $-1 =\xi_0 < \xi_1 < \ldots < \xi_N = +1$. Let $\bm U^e = [U_0^e, \ldots, U_N^e]^\top$ be the nodal values in the $e$'th element and the flux at these GLL nodes is $\bm F^e = [F_0^e, \ldots, F_N^e]^\top$. Then the semi-discrete scheme in one dimension for each component of $U$ can be written as a collocation scheme given by
\begin{equation}
\frac{\Delta x_e}{2} \od{\bm U^e}{t} + D \bm F^e + \Mm^{-1} \Rm^\top \Bm (\bm F^* - \Rm \bm F^e) = 0
\label{eq:ndg0}
\end{equation}
where $D \bm F_e$ approximates the flux derivative at the GLL points, $\Mm$ is the diagonal mass matrix containing the GLL weights $w_i$ on its diagonal, $\Bm = \textrm{diag}(-1,+1)$,  $\Rm$ is the boundary restriction operator given by $ \Rm \ \bm U^e =  [U_0^e,  \ U_N^e]^\top$,  $\bm F^{*} = [F_\emh^{*}, F_\eph^{*}]^\top$ is the numerical flux which couples the solution to the neighbouring elements $e-1$ and $e+1$, respectively, and $\Delta x_e$ is the element size. The last term containing the numerical fluxes is also called the simultaneous approximation term (SAT). In the standard DG scheme, the $D\bm F^e$ term is approximated by a differentiation matrix $\D$, leading to
\begin{equation}
\frac{\Delta x_e}{2} \od{U_i^e}{t} + \sum_{j=0}^N \D_{ij} F_j^e 
- \frac{\delta_{i0}}{w_0} (F_\emh^* - F_0^e) + \frac{\delta_{iN}}{w_N}  (F_\eph^* - F_N^e) = 0, \qquad 0 \le i \le N
\label{eq:ndg1}
\end{equation}
where $\D_{ij}$ are elements of the differentiation matrix which is explained below and $\delta_{ij}$ is the Kronecker symbol. For non-linear equations, we cannot prove any type of energy/entropy property because we cannot apply SBP on the $D\bm F^e$ terms. Following~\cite{FISHER2013353}, we can modify the scheme~\eqref{eq:ndg1} into the so called flux differencing form given by
\begin{equation}
\frac{\Delta x_e}{2} \od{U_i^e}{t} + 2\sum_{j=0}^N \D_{ij} F^{\#}_{ij}
- \frac{\delta_{i0}}{w_0} (F_\emh^* - F_0^e) + \frac{\delta_{iN}}{w_N}  (F_\eph^* - F_N^e) = 0, \qquad 0 \le i \le N
\label{eq:ndg2}
\end{equation}
where $F^{\#}_{ij} = F^{\#}(U_i^e, U_j^e)$ and  $F^{\#}$ is a symmetric flux function. The proper choice of the volumetric flux $F^{\#}$ and the surface flux $F^*$ can lead to schemes which satisfy certain energy/entropy preservation property.

The scheme~\eqref{eq:ndg2} can also be written as a split form derivative scheme~\cite{FISHER2013353}, where the derivatives $D\bm F$ of the non-linear flux in~\eqref{eq:ndg0} are computed in a special way to achieve SBP property. Define the spatial derivatives of linear ($a$), quadratic ($ab$) and cubic ($abc$) terms by the following formulae,
\begin{eqnarray*}
D_1 a &=& a_x \\
D_2(a\cdot b) &=& \half [ (ab)_x + a b_x + b a_x] \\
D_3(a\cdot b \cdot c) &=& \frac{1}{4}[ (abc)_x + a(bc)_x + b(ac)_x + c(ab)_x + bca_x + ac b_x + ab c_x]
\end{eqnarray*}
The derivatives on the right hand side are approximated at the GLL nodes using a differentiation matrix $\Dm$,
\[
a_x = \Dm \av, \qquad (a b)_x = \Dm \av \bv, \qquad (abc)_x = \Dm \av \bv \cv
\]
where
\[
\Dm_{ij} = \ell_j'(\xi_i), \qquad 0 \le i,j \le N
\]
In the above equations, the product of vectors, e.g., $\av \bv$, must be interpreted componentwise. 

We now list some standard split form schemes for the compressible Euler equations from the literature along with their flux form $F^{\#}$ followed by our proposed split form.  The flux $F^{\#}$ depends on two states $U_l, U_r$, and we use curly braces to denote the arithmetic average
\[
\avg{\phi} = \avg{\phi}_{lr} = \half (\phi_l + \phi_r)
\]
The subscript $lr$ may be suppressed in some of the following formulae for clarity.
\paragraph{Central scheme} 
\[
DF = \begin{bmatrix}
D_1(\rho u) \\
D_1 (p + \rho  u  u) \\
D_1(\rho u v) \\
D_1(\rho u  w) \\
D_1(p  u + E  u)
\end{bmatrix}, \qquad
F^{\#} = \begin{bmatrix}
\avg{\rho u} \\
\avg{p + \rho u u} \\
\avg{\rho u v} \\
\avg{\rho u w} \\
\avg{p u + E u}
\end{bmatrix}
\]
This is the canonical approach where the flux vector is differentiated and the numerical flux is just the arithmetic average of the flux at the cell boundaries. With this flux, scheme~\eqref{eq:ndg1} and~\eqref{eq:ndg2} are identical; but it does not have SBP property for non-linear equations, and it does not satisfy KEP or any other energy boundedness property.
\paragraph{Kennedy and Gruber (KG)~\cite{KENNEDY2008}}
\[
DF = \begin{bmatrix}
D_2(\rho \cdot u) \\
D_1 p + D_3(\rho \cdot u \cdot u) \\
D_3(\rho \cdot u \cdot v) \\
D_3(\rho \cdot u \cdot w) \\
D_2(p \cdot u) + D_3(\rho \cdot e \cdot u)
\end{bmatrix}, \qquad
F^{\#} = \begin{bmatrix}
\avg{\rho}\avg{u} \\
\avg{p} + \avg{\rho}\avg{u}\avg{u} \\
\avg{\rho}\avg{u}\avg{v} \\
\avg{\rho}\avg{u}\avg{w} \\
\avg{p}\avg{u} + \avg{\rho}\avg{e}\avg{u}
\end{bmatrix}, \qquad e = \frac{E}{\rho}
\]
\paragraph{Ducros~\cite{Ducros2000} }
\[
DF = \begin{bmatrix}
D_2(\rho \cdot u) \\
D_1 p + D_2(\rho u \cdot u) \\
D_2(\rho u \cdot v) \\
D_2(\rho u \cdot w) \\
D_2(p \cdot u) + D_2(E \cdot u)\end{bmatrix}, \qquad
F^{\#} = \begin{bmatrix}
\avg{\rho}\avg{u} \\
\avg{p} + \avg{\rho u}\avg{u} \\
\avg{\rho u}\avg{v} \\
\avg{\rho u}\avg{w} \\
\avg{p}\avg{u} + \avg{E}\avg{u}
\end{bmatrix}
\]
\paragraph{Shima et al. (\keep)~\cite{Shima2021}}
\[
DF = \begin{bmatrix}
D_2(\rho \cdot u) \\
D_1 p + D_3(\rho \cdot u \cdot u) \\
D_3(\rho \cdot u \cdot v) \\
D_3(\rho \cdot u \cdot w) \\
\frac{1}{\gamma-1} D_2(p \cdot u) + 
T
\end{bmatrix}, \qquad
F^{\#} = \begin{bmatrix}
\avg{\rho}\avg{u} \\
\avg{p} + \avg{\rho}\avg{u}\avg{u} \\
\avg{\rho}\avg{u}\avg{v} \\
\avg{\rho}\avg{u}\avg{w} \\
\frac{1}{\gamma-1} \avg{p}\avg{u} + \half \avg{\rho}(u_l u_r + v_l v_r + w_l w_r)\avg{u} +  \half (p_l u_r + p_r u_l)
\end{bmatrix}
\]
where $T$ is the split form for the kinetic energy and the $pu$ terms.
\paragraph{Modified KEP (mKEP)} Finally, our proposed modified KEP scheme is given by
\[
DF = \begin{bmatrix}
D_2(\rho \cdot u) \\
D_1 p + D_3(\rho \cdot u \cdot u) \\
D_3(\rho \cdot u \cdot v) \\
D_3(\rho \cdot u \cdot w) \\
\frac{\gamma}{\gamma-1} D_2(p \cdot u) +  D_3(\rho \cdot k \cdot u)
\end{bmatrix}, \qquad
F^{\#} = \begin{bmatrix}
\avg{\rho}\avg{u} \\
\avg{p} + \avg{\rho}\avg{u}\avg{u} \\
\avg{\rho}\avg{u}\avg{v} \\
\avg{\rho}\avg{u}\avg{w} \\
\frac{\gamma}{\gamma-1} \avg{p}\avg{u} + \avg{\rho}\avg{k}\avg{u}
\end{bmatrix}, \qquad k = \half (u^2 + v^2 + w^2)
\]
This is similar to the KG flux~\cite{KENNEDY2008} and Jameson flux~\cite{Jameson2008} in all the components except the energy equation, where we use only primitive variables and treat the flux in quadratic and cubic form.

The KG, \keep and mKEP schemes are kinetic energy preserving, but the central and Ducros are not.
\subsection{SBP property} 
The differentiation matrix $\D$ satisfies the summation-by-parts (SBP) property given by
\begin{equation*}
\Mm \Dm + \Dm^\top \Mm = \textrm{diag}(-1,0,\ldots,0,+1), \qquad \Mm = \textrm{diag}(w_0, w_1, \ldots, w_N)
\end{equation*}
which can be written in component form as
\begin{equation}
w_i \D_{ij} + w_j \D_{ji} = - \delta_{i0} \delta_{j0} + \delta_{iN} \delta_{jN}, \qquad 0 \le i,j \le N
\label{eq:sbp}
\end{equation}
The SBP property implies the following result
\[
\bm f^\top \Mm \Dm \bm g + \bm f^\top \Dm^\top \Mm \bm g = f_N g_N  - f_0 g_0
\]
which is the discrete analogue of the integration by parts formula
\[
\int_{-1}^{+1} (f g_\xi + g f_\xi) \ud \xi =  f(+1) g(+1) - f(-1) g(-1)
\]
Let us also note the following SBP result for later use. If $f_{ij}$ is some symmetric quantity, i.e., $f_{ij} = f_{ji}$, then it is easy to show that
\begin{equation}
\sum_{i=0}^N \sum_{j=0}^N w_i \D_{ij} f_{ij} = -f_{00} + f_{NN}
\label{eq:sbp2}
\end{equation}
\section{One dimensional density wave}
\label{sec:dwave}
Consider a simple solution to the 1-D Euler equations given by~\cite{gassner2020stability} 
\[
u(x,t) = V = \textrm{constant}, \qquad p(x,t) = P = \textrm{constant}
\]
and where the density $\rho(x,t) = r(x,t)$ solves the linear advection equation
\[
\ro_t + V \ro_x = 0
\]
with periodic boundary conditions. We first check if the split form schemes are able to maintain the constancy of velocity and pressure when applied to solve this problem, by writing the schemes in split derivative form.
\subsection{Behaviour of the KG scheme}
Consider the initial condition corresponding to the density wave solution. The split flux derivatives at the initial time are given by
\[
D_2(\rho \cdot u) = V \Dm \rhov, \qquad D_1 p = 0, \qquad D_2(\rho \cdot u \cdot u) = V^2 \Dm \rhov
\]
\[
D_2(p \cdot u) = 0, \qquad D_3(\rho\cdot e \cdot u) = \frac{V^3}{2} \Dm \rhov + \frac{P V}{2} \rhov \Dm (1/\rhov) + \frac{P V}{2} (1/\rhov) \ \Dm \rhov
\]
so that
\begin{eqnarray*}
\od{\rhov}{t}|_{t=0} &=& - V \Dm \rhov + \SAT_\rho \\
\od{\mv}{t}|_{t=0} &=& - V^2 \Dm \rhov + \SAT_{\rho u} \\
\od{\Ev}{t}|_{t=0} &=& - \frac{V^3}{2} \Dm \rhov - \frac{P V}{2}\left[  \rhov \Dm (1/\rhov) + (1/\rhov) \Dm \rhov \right] + \SAT_E,
\end{eqnarray*}
where $\SAT_\rho, \SAT_{\rho u}, \SAT_E$ denote the $\SAT$ terms for the discretization of $\rho, \rho u, E$ respectively. At the continuous level $\rho (1/\rho)_x + (1/\rho) \rho_x = 0$, but at the discrete level $\rhov \Dm (1/\rhov) + (1/\rhov) \Dm \rhov \ne 0$. The velocity and pressure equation can be obtained from
\[
\od{\uv}{t}|_{t=0} = \frac{1}{\rhov} \left[ \od{\mv}{t} - V \od{\rhov}{t} \right]_{t=0} = 0, \qquad \od{\pv}{t}|_{t=0} = (\gamma-1) \left[ \od{\Ev}{t} - V \od{\mv}{t} + \half V^2 \od{\rhov}{t} \right]_{t=0} \ne 0
\]
The scheme cannot maintain constancy of pressure, and consequently also the velocity, and it is unstable for this problem in numerical computations.
\subsection{Behaviour of the modified KEP scheme}
At the initial condition for the density wave solution, we have similar relations as in previous sub-section and
\[
D_3(\rho\cdot k \cdot u) = \frac{V^3}{2} \Dm \rhov
\]
so that
\[
\od{\rhov}{t}|_{t=0} = - V \Dm \rhov + \SAT_\rho, \qquad 
\od{\mv}{t}|_{t=0} = - V^2 \Dm \rhov + V \ \SAT_\rho, \qquad
\od{\Ev}{t}|_{t=0} = - \frac{V^3}{2} \Dm \rhov  + \frac{V^3}{2} \SAT_\rho
\]
which implies that
\[
\od{\uv}{t}|_{t=0} = 0, \qquad \od{\pv}{t}|_{t=0} = 0
\]
The scheme maintains constancy of velocity and pressure. 
\subsection{Behaviour of some other schemes}
The central, Ducros and \keep schemes are also able to maintain constancy of pressure and velocity. However, other schemes like Kennedy and Gruber (KG), Jameson~\cite{Jameson2008}, Morinishi\cite{MORINISHI2010276}, Kuya et al.~\cite{Kuya2018}, etc., which make use of enthalpy or specific internal energy or specific total energy in the energy flux and thereby cause a coupling of density and pressure, do not maintain constancy of pressure and thereby the velocity. Such schemes are also found to be unstable in numerical computations for this test problem.
\subsection{Linearized equations}
We now consider a small perturbation around the density wave solution and write the resulting density, velocity and pressure as $\ro+\rho, V+u, P+p$. Here $(\rho,u,p)$ are small perturbations which are governed by the linearized Euler equations given by
\begin{eqnarray}
\nonumber
\rho_t + V \rho_x + (\ro u)_x &=& 0 \\
\label{eq:lineul}
u_t + V u_x + \frac{1}{\ro} p_x &=& 0 \\
\nonumber
p_t + V p_x + \gamma P u_x &=& 0
\end{eqnarray}
These equations admit a reduced energy conservation equation of the form
\begin{equation}
\pen_t + (\pen V + pu)_x = 0, \qquad \pen = \frac{p^2}{2\gamma P} + \half \ro u^2
\label{eq:peneq}
\end{equation}
so that we have the conservation of the total energy in the domain
\begin{equation}
\int\pen\ud x = \textrm{constant}
\label{eq:pencon}
\end{equation}
under periodic boundary condition. The energy $\pen$ does not involve the density perturbations and hence we refer to it as a reduced energy. We are not able to derive an energy conservation equation that includes all three variables. See Appendix \ref{sec:sym} for an alternate approach based on symmetrization where an energy equation is established but this energy is not strictly bounded with respect to time. The {\em energy} in the density perturbations is given by the equation
\[
(\rho^2/2)_t + V (\rho^2/2)_x + \rho (\ro u)_x = 0
\]
so that
\begin{equation}
\od{}{t}\int \half \rho^2 \ud x = - \int \rho (\ro u)_x \ud x
\label{eq:rhoperttot}
\end{equation}
The convective term does not change the total density perturbations in the domain. The density perturbations however can grow with time if the right hand side is positive and we do not have strict control on the density perturbations.
\subsection{Linearization of modified KEP scheme}
Let $\ro,V,P$ be solution of the modified KEP scheme with constant velocity and pressure. Note that we are analyzing the semi-discrete scheme where time is still continuous. So $r$ satisfies the ODE
\begin{equation}
\frac{\Delta x_e}{2}
\od{r_i^e}{t} + 2 \sum_{j=0}^N \D_{ij} V \avg{r}_{ij} - \frac{\delta_{i0}}{w_0} V(\avg{\ro}_\emh - \ro_0^e) + \frac{\delta_{iN}}{w_N} V  (\avg{r}_\eph - \ro_N^e) = 0
\label{eq:roode}
\end{equation}
We consider a linearization around this solution. The linearized density equation is
\begin{equation}
\begin{aligned}
\frac{\Delta x_e}{2}
\od{\rho_i^e}{t} + 2 \sum_{j=0}^N \D_{ij} [V\avg{\rho}_{ij} + \avg{\ro}_{ij} \avg{u}_{ij}] - & \frac{\delta_{i0}}{w_0} [f_\emh^{*\rho} - (\ro u + \rho V)_0^e] + \frac{\delta_{iN}}{w_N}  [f_\eph^{*\rho} - (\ro u + \rho V)_N^e] = 0 \\
f^{*\rho} =& \ \avg{\ro} \avg{u} + V \avg{\rho}
\end{aligned}
\label{eq:rlinode}
\end{equation}
To first order, the momentum time derivative is given by
\[
\od{}{t}(\ro_i^e + \rho_i^e)(V + u_i^e) \approx V \od{\ro_i^e}{t} + V \od{\rho_i^e}{t} + \od{}{t}(\ro_i^e u_i^e)
\]
and using~\eqref{eq:roode}, the linearized momentum equation can be written as
\begin{equation}
\begin{aligned}
\frac{\Delta x_e}{2}
\od{}{t}(\ro_i^e u_i^e) + 2 \sum_{j=0}^N \D_{ij}[ \avg{p}_{ij} + V \avg{\ro}_{ij} \avg{u}_{ij} ] -& \frac{\delta_{i0}}{w_0} [f_\emh^{*m} - (p + \ro V u)_0^e] + \frac{\delta_{iN}}{w_N}  [f_\eph^{*m} - (p + \ro V u)_N^e] = 0 \\
f^{*m} =& \ \avg{p} + V \avg{\ro} \avg{u}
\end{aligned}
\label{eq:mlinode}
\end{equation}
To first order, the energy time derivative is
\[
\od{E_i^e}{t} \approx \frac{V^2}{2} \od{\ro_i^e}{t} + \frac{V^2}{2} \od{\rho_i^e}{t} + V \od{}{t}(\ro_i^e u_i^e) + \frac{1}{\gamma-1} \od{p_i^e}{t}
\]
Then using~\eqref{eq:rlinode},~\eqref{eq:mlinode}, the linearized pressure equation is obtained as
\begin{equation}
\begin{aligned}
\frac{\Delta x_e}{2}
\od{p_i^e}{t} + 2 \sum_{j=0}^N \D_{ij}[ V \avg{p}_{ij} + \gamma P \avg{u}_{ij} ] -& \frac{\delta_{i0}}{w_0} [f_\emh^{*p} - (p V + \gamma P u)_0^e] + \frac{\delta_{iN}}{w_N}  [f_\eph^{*p} - (p V + \gamma P u)_N^e] = 0 \\
f^{*p} = & \ V \avg{p} + \gamma P \avg{u}
\end{aligned}
\label{eq:linpre}
\end{equation}
Then the perturbation energy $\pen$ defined in~\eqref{eq:peneq} satisfies
\[
\od{\pen}{t} = \frac{p}{\gamma P} \od{p}{t} + \od{}{t}  \half \ro u^2 = \frac{p}{\gamma P} \od{p}{t} + u \od{}{t}(\ro u) - \frac{u^2}{2} \od{\ro}{t}
\]
Using the above identity and~\eqref{eq:roode},~\eqref{eq:mlinode},~\eqref{eq:linpre}, we can show that the energy at a solution node satisfies
\[
\frac{\Delta x_e}{2}
\od{\pen_i^e}{t} + 2 \sum_{j=0}^N \D_{ij} \left[ \frac{V}{\gamma P} p_i \avg{p}_{ij} + p_i \avg{u}_{ij} + u_i \avg{p}_{ij} + V u_i \avg{r}_{ij} \avg{u}_{ij} - \half V u_i^2 \avg{\ro}_{ij} \right] - A_i + B_i = 0
\]
where $A_i,B_i$ are the interface terms at left and right faces of the element. This equation can be re-written as
\[
\frac{\Delta x_e}{2}
\od{\pen_i^e}{t} + \sum_{j=0}^N \D_{ij} \left[ \frac{V}{\gamma P} p_i p_j + p_i u_j + u_i p_j +  V u_i u_j \avg{r}_{ij} \right]  - A_i + B_i = 0
\]
We multiply by the quadrature weight $w_i$ and sum over all $i$; using the symmetry of the quantity inside the square brackets and the SBP property~\eqref{eq:sbp}, we obtain
\begin{align*}
\frac{\Delta x_e}{2}
\od{}{t} \sum_i w_i \pen_i^e - \half \left[ \frac{V}{\gamma P} (p_0^e)^2 + 2 p_0^e u_0^e +  V (u_0^e)^2 \ro_0^e \right] & + \half \left[ \frac{V}{\gamma P} (p_N^e)^2 + 2 p_N^e u_N^e +  V (u_N^e)^2 \ro_N^e \right] \\
& - \sum_i w_i A_i + \sum_i w_i B_i = 0
\end{align*}
The boundary term can be written as
\begin{eqnarray*}
\sum_i w_i A_i &=& \frac{p_0^e}{\gamma P} f_\emh^{*p} + u_0^e f_\emh^{*m} - \half (u_0^e)^2 V \avg{\ro}_\emh - \left[ \frac{p_0^e}{\gamma P} (p_0 V + \gamma P u_0) + u_0^e (p_0 + \ro_0 V u_0) - \half (u_0^e)^2 V \ro_0^e \right]\\
&=& \left[ \frac{V}{\gamma P} p_0^e \avg{p}_\emh + p_0^e \avg{u}_\emh + u_0^e \avg{p}_\emh + \half V u_N^{e-1} u_0^e \avg{\ro}_\emh \right]\\
&& - \left[ \frac{V}{\gamma P} (p_0^e)^2 + 2 p_0^e u_0^e +  \half V (u_0^e)^2 \ro_0^e \right]
\end{eqnarray*}
with a similar equation for the other interface term, so that the energy equation becomes
\begin{align*}
\frac{\Delta x_e}{2}
\od{}{t} \sum_i w_i \pen_i^e - & \half \left[ \frac{V}{\gamma P} p_N^{e-1} p_0^e +  p_0^e u_N^{e-1} +  p_N^{e-1} u_0^e + V u_N^{e-1} u_0^e \avg{\ro}_\emh \right] \\
+ &  \half \left[ \frac{V}{\gamma P} p_N^{e} p_0^{e+1} +  p_0^{e+1} u_N^{e} +  p_N^{e} u_0^{e+1} + V u_N^{e} u_0^{e+1} \avg{\ro}_\eph \right] = 0
\end{align*}
This has a divergence form; when we add the equations from all the elements, there is a telescopic collapse of the terms inside the square brackets leading to
\[
\od{}{t} \sum_e \sum_i w_i \pen_i^e = 0
\]
The total perturbation energy is conserved which is consistent with the continuous case of equation~\eqref{eq:pencon}.

From \eqref{eq:rlinode}, we obtain an equation for the {\em energy} in the density perturbations, which is given by
\begin{equation*}
\frac{\Delta x_e}{2}
\sum_i w_i \rho_i^e \od{\rho_i^e}{t} + \sum_i \sum_{j=0}^N w_i \D_{ij} [V \rho_i \rho_j  + 2 \rho_i \avg{\ro}_{ij} \avg{u}_{ij}] - \rho_0^e [f_\emh^{*\rho} - (\ro u + \rho V)_0^e] + \rho_N^e [f_\eph^{*\rho} - (\ro u + \rho V)_N^e] = 0
\end{equation*}
Using SBP property~\eqref{eq:sbp2}, we get
\[
\sum_i \sum_j w_i \D_{ij} \rho_i^e \rho_j^e = -\half (\rho_0^e)^2 + \half (\rho_N^e)^2
\]
and adding the equation from all elements, it is easy to show that the convective terms involving $V$ cancel out just as in the continuous analysis.
\paragraph{Remark.}
The \keep flux also leads to the same linearized equations as the mKEP flux. It behaves similarly to the mKEP flux when numerically computing small perturbations as shown in the results.
\subsection{Linearization of central scheme}
The base flow density $r$ evolves according to~\eqref{eq:roode}.  The linearized density equation is
\begin{equation}
\frac{\Delta x_e}{2}
\od{\rho_i^e}{t} + 2 \sum_{j=0}^N \D_{ij} [V\avg{\rho}_{ij} + \avg{\ro u}_{ij}] - \frac{\delta_{i0}}{w_0} [f_\emh^{*\rho} - (\ro u + \rho V)_0^e] + \frac{\delta_{iN}}{w_N}  [f_\eph^{*\rho} - (\ro u + \rho V)_N^e] = 0
\label{eq:rlinodec}
\end{equation}
\[
f^{*\rho} = \avg{\ro u} + V \avg{\rho}
\]
and the momentum equation is
\begin{equation}
\frac{\Delta x_e}{2}
\od{}{t}(\ro_i^e u_i^e) + 2 \sum_{j=0}^N \D_{ij}[ \avg{p}_{ij} + V \avg{\ro u}_{ij} ] - \frac{\delta_{i0}}{w_0} [f_\emh^{*m} - (p + \ro V u)_0^e] + \frac{\delta_{iN}}{w_N}  [f_\eph^{*m} - (p + \ro V u)_N^e] = 0
\label{eq:mlinodec}
\end{equation}
\[
f^{*m} = \avg{p} + V \avg{\ro u}
\]
The pressure equation is same as in the case of mKEP scheme and is given by~\eqref{eq:linpre}. We then obtain the following energy equation
\[
\frac{\Delta x_e}{2}
\od{\pen_i^e}{t} + \frac{1}{4}V \sum_{j=0}^N \D_{ij} u_i (r_j - r_i)(u_j - u_i) + \sum_{j=0}^N \D_{ij} \left[ \frac{V}{\gamma P} p_i p_j + p_i u_j + u_i p_j +  V u_i u_j \avg{r}_{ij} \right]  - A_i + B_i = 0
\]
The third term is same as in the previous section and we can apply SBP on it. The second term cannot be simplified using SBP property and leads to a volumetric term with non-determinate sign in the energy equation. This is more clearly shown in the finite difference context in Appendix~\ref{sec:fd}. Thus the central scheme is not consistent with the energy equation~\eqref{eq:pencon}. A similar result is obtained for the case of Ducros flux.
\paragraph{Remark.} Appendix \ref{sec:fd} provides linearized schemes in the finite difference context where it is easier to do the analysis.
\section{Numerical tests}
In this section,
we first present numerical results for the compressible Euler equations for three different test cases. 
We use the open source code Trixi.jl~\cite{trixi}.  Time-stepping is performed with the fourth-order, five-stage, low-storage Runge-Kutta method of \cite{carpenter1994} and with CFL=0.2. 
The interface flux $F^*$ is taken to be same as the volumetric symmetric flux $F^\#$, i.e., $F^*_\eph = F^\#(U_N^e, U_0^{e+1})$, since our purpose here is to investigate the properties of the scheme when there is no added numerical dissipation.  
All the meshes used 
for compressible Euler computations are comprised of uniform quadrilateral  or hexahedral elements.
Next, we present implicit large eddy simulations (ILES) for 
the compressible Navier Stokes equation
for two standard test cases. 
These simulations are performed using the solver documented in~\cite{SINGH2020, Singh2021}.
\subsection{2-D density wave}
Consider the initial condition given by
\[
\rho = 1 + 0.98 \sin(2\pi(x+y)), \qquad u = 0.1, \qquad v = 0.2, \qquad p = 20.0
\]
in the domain $[-1,+1]^2$ with periodic boundary conditions. The exact solution is a translation of the density profile at constant velocity. This is an exact, smooth solution for the compressible Euler equations. However, as shown in \cite{gassner2020stability, ranocha2020preventing}, many popular schemes, including entropy conserving schemes blow-up for this solution due to loss of positivity of density. We run the test until time $T=100$ and monitor the time at which blow-up happens. Table~\ref{tab:den1} shows the time to blow up using mesh of $4 \times 4$ cells and $N=3,4,5$,  while Table~\ref{tab:den2} shows results for $8 \times 8$ mesh.  Note that unstable cases were repeated with CFL=0.05 but  this did not change the results.  All schemes show nearly identical time to blow-up except KG which blows up much earlier, even with increased resolution. These can be attributed to non-preservation of constant velocity and pressure condition by the scheme and the presence of unstable eigenvalues for KG flux as shown in~\cite{gassner2020stability}.


\begin{table}[!h]
\begin{center}
\begin{tabular}{|c|c|c|c|c|c|}
\hline
N  & Central  & Ducros & KG   &  \keep & mKEP \\
  \hline
3 & 0.51 & 0.51   & 0.13 & 0.51 & 0.51\\
\hline
4 & 0.49 & 0.49   & 0.07 & 0.49 & 0.49\\
\hline
5 & -   & -     & 0.08 & -  & -\\
\hline
\end{tabular}
\end{center}
\caption{Time to blow-up for density wave test for $4 \times 4$ mesh. Tests were run until time $T=100$.}
\label{tab:den1}
\end{table}

\begin{table}[!h]
\begin{center}
\begin{tabular}{|c|c|c|c|c|c|}
\hline
N  & Central   & Ducros & KG   & \keep & mKEP\\
  \hline
3 & - & -   & 0.09 & - & -\\
\hline
4 & - & -   & 0.08 & - & - \\
\hline
5 & -   & -     & 0.08 & -  & -\\
\hline
\end{tabular}
\end{center}
\caption{Time to blow-up for density wave test for $8 \times 8$ mesh. Tests were run until time $T=100$.}
\label{tab:den2}
\end{table}

Next we add a perturbation to the solution, so that the initial velocity is given by
\\
\[
u = 0.1 + A[ \sin(2\pi x) + \sin(2\pi y)], \qquad v = 0.2 + A [ \cos(2\pi x) + \cos(2\pi y)]
\]
where $A$ is the amplitude of the perturbation. Tables \ref{tab:den_per1}, \ref{tab:den_per2} and \ref{tab:den_per3} show time to blow-up for  $4 \ \times \ 4$ mesh with $N=5$,  $ 8 \ \times \ 8$ mesh with $N=3$, and $ 8 \ \times \ 8$ mesh with $N=4$, respectively. Stronger perturbations with $A=10^{-3}$ result in small time  to blow-up and maybe indicative of nonlinear instabilities. With $A=10^{-4}$, we can see simulations run for much longer times until blow-up. For this amplitude, mKEP and \keep always run as long as or longer than Central and Ducros  schemes.  This is consistent with our linear stability analysis as mKEP and \keep satisfy the reduced energy conservation equation and are also KEP.

\begin{table}[!h]
\begin{center}
\begin{tabular}{|c|c|c|c|c|c|}
\hline
A & Central   & Ducros & KG   & \keep & mKEP\\
  \hline
1E-3 & 0.80 & 7.48   & 0.08 & 5.79 & 5.79\\
\hline
1E-4 & 40.22 & 40.22   & 0.08 & 49.17 & 49.17\\
\hline
1E-5 & -   & -     & 0.08 & -  & -\\
\hline
\end{tabular}
\end{center}
\caption{Time to blow-up for density wave test with perturbation for $4 \times 4$ mesh, N=5. Tests were run until time $T=100$}
\label{tab:den_per1}
\end{table}

\begin{table}[!h]
\begin{center}
\begin{tabular}{|c|c|c|c|c|c|}
\hline
A  & Central   & Ducros & KG   & \keep & mKEP \\
  \hline
1E-3 & 3.90 & 3.90   & 0.07 & 2.47 & 2.47\\
\hline
1E-4 & 20.83 & 20.83   & 0.10 & 20.84 & 20.84\\
\hline
1E-5 & -   & -     & 0.10 & -  & -\\
\hline
\end{tabular}
\end{center}
\caption{Time to blow-up for density wave test with perturbation for $8 \times 8$ mesh, N=3. Tests were run until time $T=100$.}
\label{tab:den_per2}
\end{table}

\begin{table}[!h]
\begin{center}
\begin{tabular}{|c|c|c|c|c|c|}
\hline
A  & Central   & Ducros & KG   & \keep & mKEP\\
  \hline
1E-3 & 2.50 & 9.13   & 0.08 & 9.88 & 9.88\\
\hline
1E-4 & 77.47 & 84.56   & 0.09 & - & -\\
\hline
1E-5 & -  & -     & 0.09 & -  & -\\
\hline
\end{tabular}
\end{center}
\caption{Time to blow-up for density wave test with perturbation for $8 \times 8$ mesh, N=4. Tests were run until time $T=100$.}
\label{tab:den_per3}
\end{table}
\subsection{2-D isentropic vortex}

The isentropic vortex is another exact solution of the compressible Euler equations, and is a popular test case for evaluating high order methods.  Various initial conditions are available in literature \cite{Spiegel2016} for this problem. We use the initial condition given by
\[ 
\rho = \left[ 1 - \frac{\beta^2 (\gamma - 1)}{8 \gamma \pi^2} \exp (1 -
   r^2) \right]^{\frac{1}{\gamma - 1}}, \qquad u = M \cos \alpha -
   \frac{\beta (y - y_c)}{2 \pi} \exp \left( \frac{1 - r^2}{2} \right) 
\]
\[ v = M \sin \alpha + \frac{\beta (x - x_c)}{2 \pi} \exp \left( \frac{1 -
   r^2}{2} \right), \qquad r^2 = (x - x_c)^2 + (y - y_c)^2 
\]
and the pressure is given by $p = \rho^{\gamma}$. We choose the parameters $\beta = 5$, $M = 0.5$, $\alpha = 45^o$, $(x_c, y_c) = (0, 0)$ and the domain is taken to be $[- 10, 10] \times [- 10, 10]$.  We run the computations upto a time of $t = 2 \sqrt{2} \cdot 20 / M = 113.13708499$ when the vortex has crossed the domain two times in the diagonal direction. To measure the conservation of total kinetic energy and entropy, we plot the relative change in these quantities defined as
$$
ke(t) = \frac{KE(t) - KE(0)}{KE(0)}, \qquad
en(t) = \frac{EN(t) - EN(0)}{|EN(0)|},
$$
where $KE(t), EN(t)$ refers to total kinetic energy and entropy in the spatial domain at time $t$. The entropy refers to the mathematical entropy which is a convex function that should be non-increasing with time. Note that we divide by $|EN(0)|$ since the total entropy can be negative; the non-increasing property of entropy requires that $en(t)$ be a non-increasing function of time.

\begin{figure}
  \begin{center}
    \includegraphics[width=0.98\textwidth]{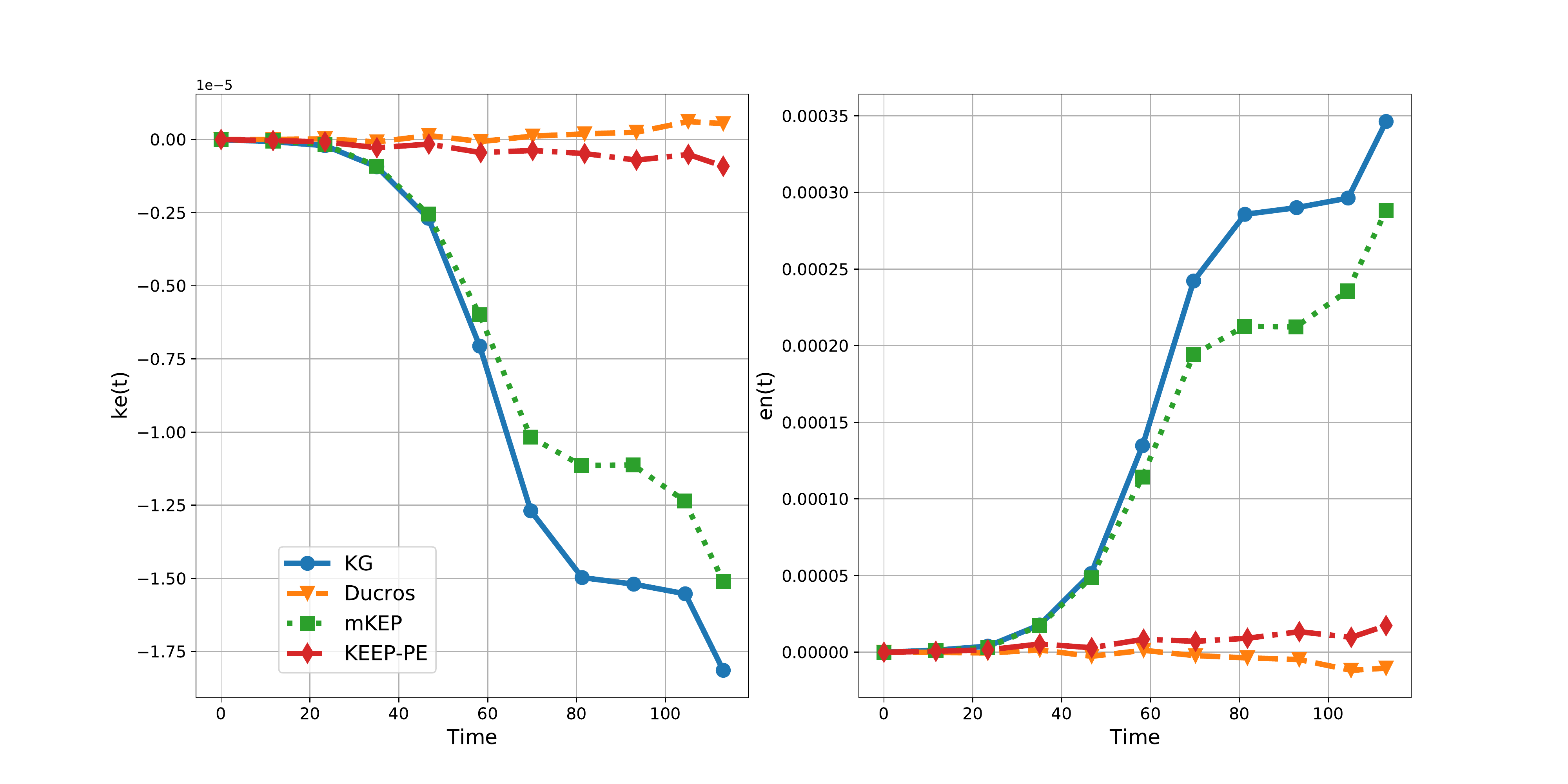}
  \caption{Comparison of kinetic energy and entropy variation with time for isentropic vortex on $32 \times 32$ mesh with $N=3$.}
  \label{f:iv_n_3}
  \end{center}
\end{figure}
\begin{figure}
  \begin{center}
    \includegraphics[width=0.98\textwidth]{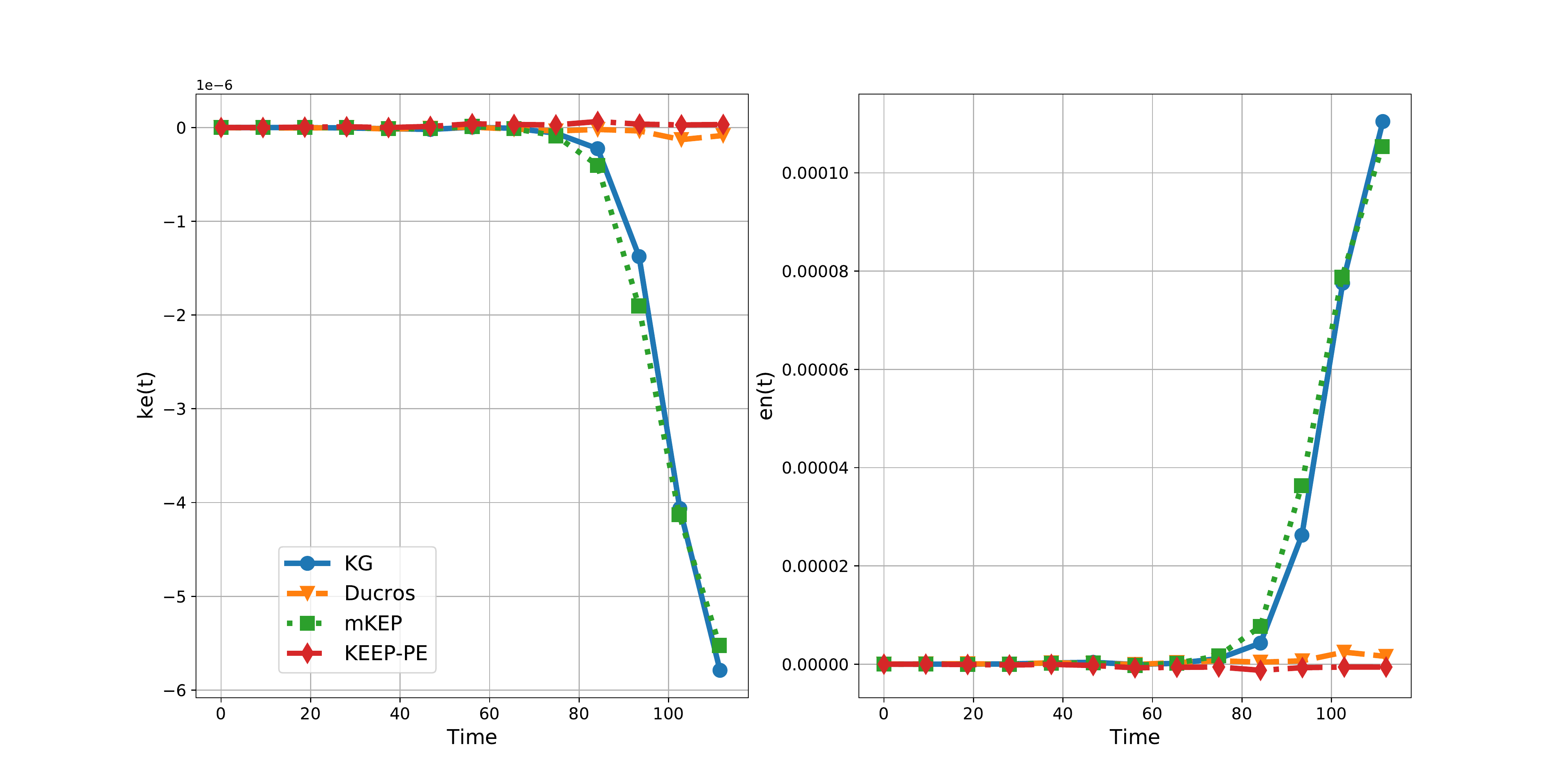}
  \caption{Comparison of kinetic energy and entropy variation with time for isentropic vortex on $32 \times 32$ mesh with $N=4$.}
  \label{f:iv_n_4}
  \end{center}
\end{figure}

The simulations were performed using a $32 \ \times \ 32$  uniform quadrilateral mesh with polynomial orders of $N=3, 4$.  Figure~\ref{f:iv_n_3} and \ref{f:iv_n_4} show the evolution of  normalized kinetic energy and entropy for the various split forms. The central scheme was not stable for these simulations and is not shown in the figures.  We observe that KG and mKEP behave very similarly;  they are both kinetic energy dissipative and mildly entropy unstable for all simulation cases.  In contrast, both Ducros and \keep schemes are nearly conservative for both  kinetic energy and entropy. We observe slight kinetic energy growth for Ducros and slight entropy growth for \keep. For $N=4$, the conservation violations are even smaller. 
\subsection{3-D inviscid Taylor-Green vortex}
\label{sec:tgv1}
The density wave and isentropic vortex test cases were 2D cases with smooth solutions, so that turbulence effects on flow stability do not exist in those problems. We use the 3D Taylor-Green vortex case to test the numerical properties of the scheme in the presence of turbulence. Here we use the inviscid version which is useful in studying the  stability of the discrete formulation~\cite{Gassner2016, Shima2021} for the Euler equations.  The domain is a cube of length $2\pi$ and the initial conditions are given by
%
\[
\rho = \rho_0, \qquad u = +U_0 \sin(x) \cos(y) \cos(z), \qquad
v = -U_0 \cos(x) \sin(y) \cos(z), \qquad
w = 0
\]
\[
p = P_0 + 
\frac{\rho_0 U_0^2}{16}(\cos(2x)\cos(2z) + 2\cos(2x) + 2\cos(2y) + \cos(2y)\cos(2z) ), \qquad
P_0 = \frac{\rho_0}{M^2 \gamma}
\]
with $M$ being the Mach number and $\rho_0 = 1, U_0 = 1, M = 0.1$. We run the simulations on a $16 \times 16 \times 16$ uniform hexahedral mesh with $N=3$.
\begin{figure}
  \begin{center}
    \includegraphics[width=0.99\textwidth]{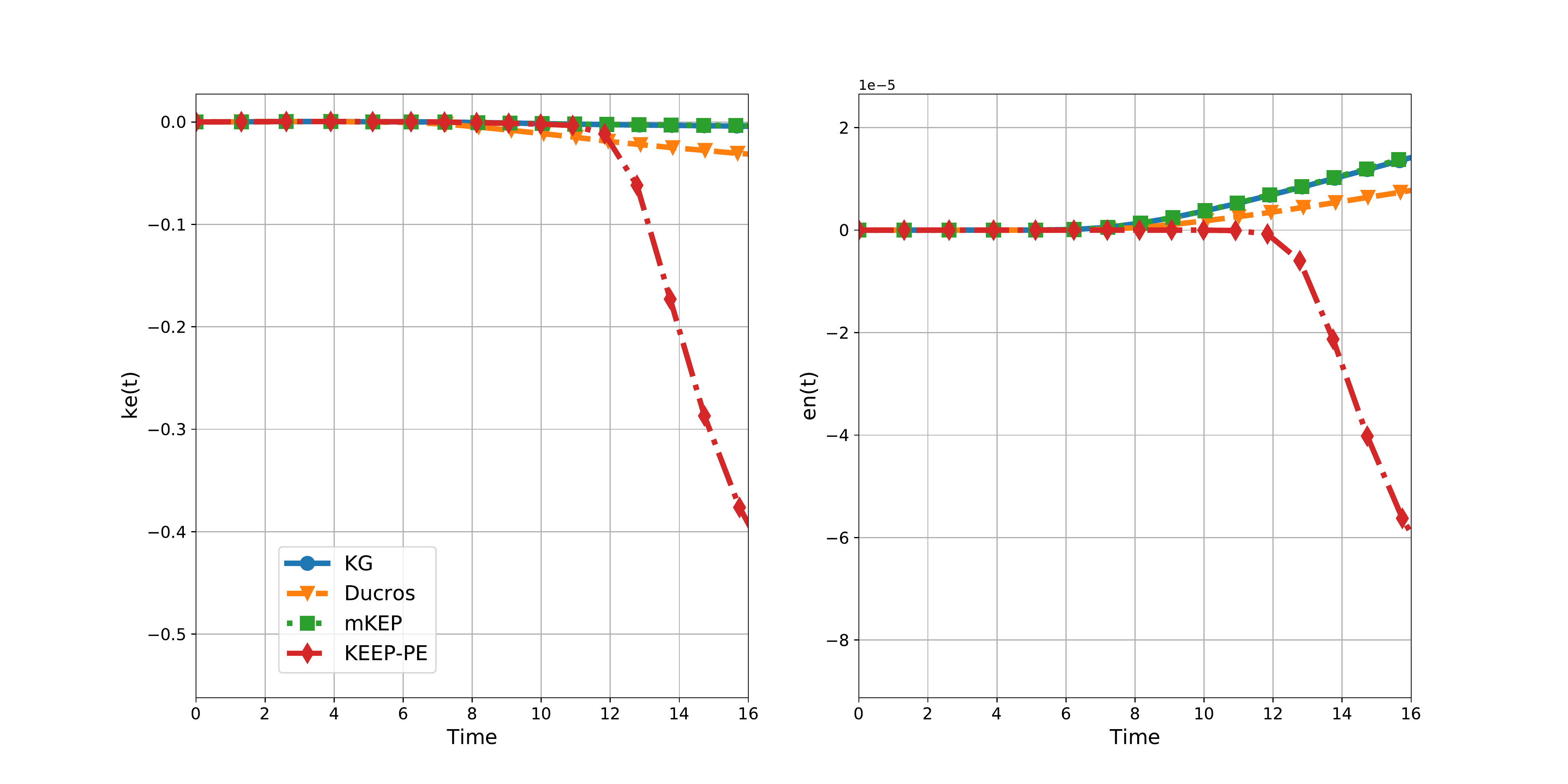}
  \caption{Comparison of kinetic energy and entropy variation with time for inviscid TGV. 
  Simulations were run on a $16 \times 16 \times 16$ uniform hexahedral mesh with $N=3$
  till T=20.
  Note the scale for entropy.}
  \label{f:tgv}
  \end{center}
\end{figure}
\begin{figure}
  \begin{center}
    \includegraphics[width=0.99\textwidth]{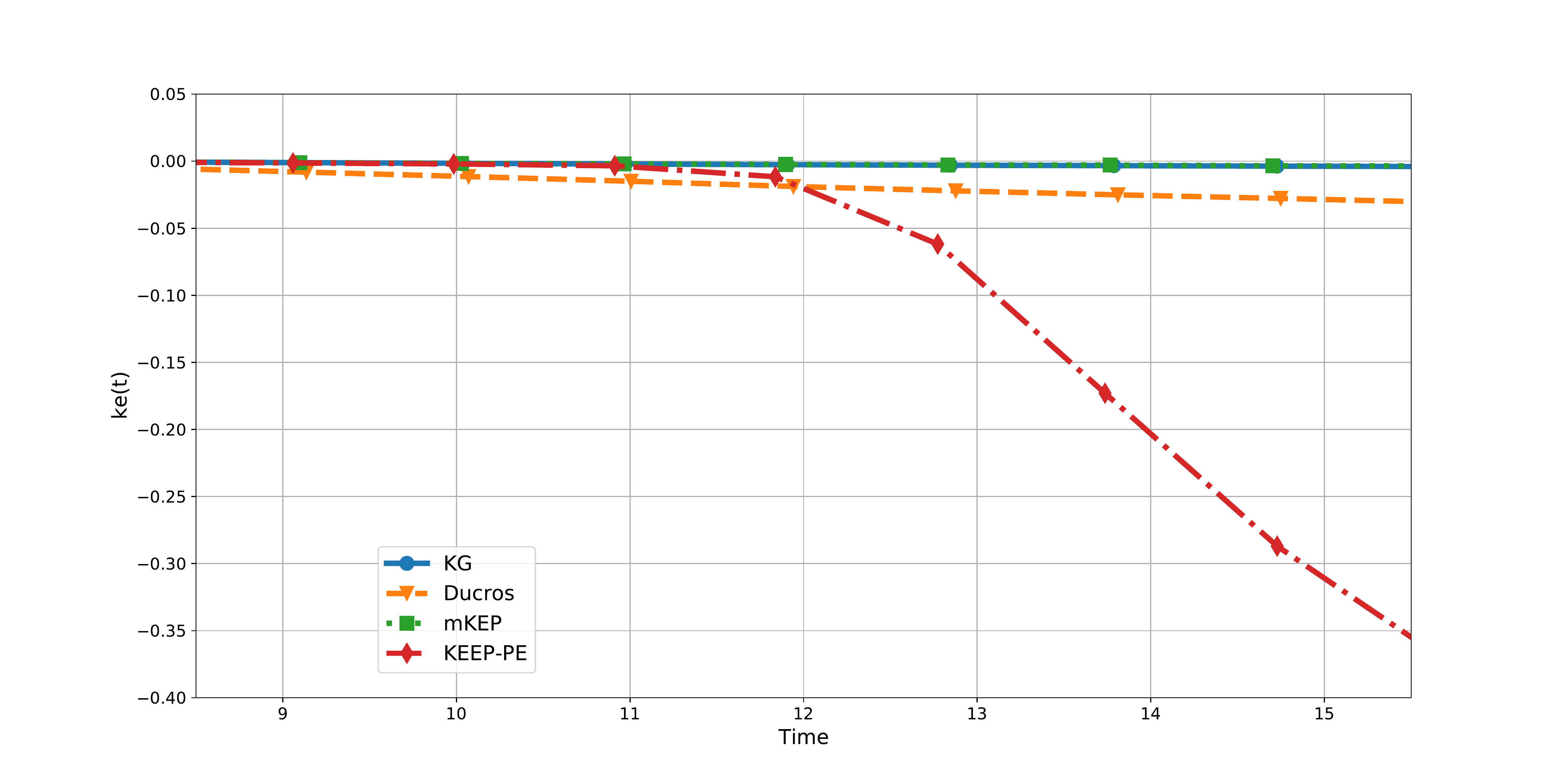}
  \caption{Zoomed in view of kinetic energy variation for inviscid TGV. Note that KG and mKEP are nearly kinetic energy conservative compared to the other schemes.}
  \label{f:tgv_zoom}
  \end{center}
\end{figure}

Figure~\ref{f:tgv} shows the evolution of normalized kinetic energy and entropy for the different split forms. As with the isentropic vortex case, central scheme was not stable. All split forms nearly conserve kinetic energy and entropy until about $T=6$, when the flow starts transitioning and becomes turbulent. The mKEP and KG split forms show nearly identical results and have very low kinetic energy dissipation. The Ducros scheme has slightly higher dissipation levels, while the \keep scheme has very high dissipation of kinetic energy after the onset  of turbulence.  Figure~\ref{f:tgv_zoom} shows a zoomed in view of kinetic energy evolution  where the advantage of the mKEP and KG schemes over the other two schemes become even more  apparent.

In the case of total entropy, KG and mKEP schemes are again nearly identical and are mildly entropy unstable; Ducros is also entropy unstable but less than mKEP, and \keep is the only scheme that is entropy stable.  However, note again that the entropy instability is mild for mKEP and \keep schemes.  This is in contrast to the kinetic energy dissipation which is  substantial for the \keep split form.

\subsection{Implicit LES results}

High order DG methods and especially their split form versions are suitable for implicit large eddy simulation (ILES) of turbulent flows because of the high wave number dissipation properties of the method,  which works as an LES model \cite{Beck2014, Flad2017, Winters2018, SINGH2020}. To demonstrate the capabilities of the new split form, we present results for two benchmark test cases. In these simulations, we add Lax-Friedrich dissipation to the surface flux between the elements. 

\subsubsection{3-D viscous Taylor-Green vortex}

\begin{figure}
  \begin{center}
    \includegraphics[width=1.05\textwidth]{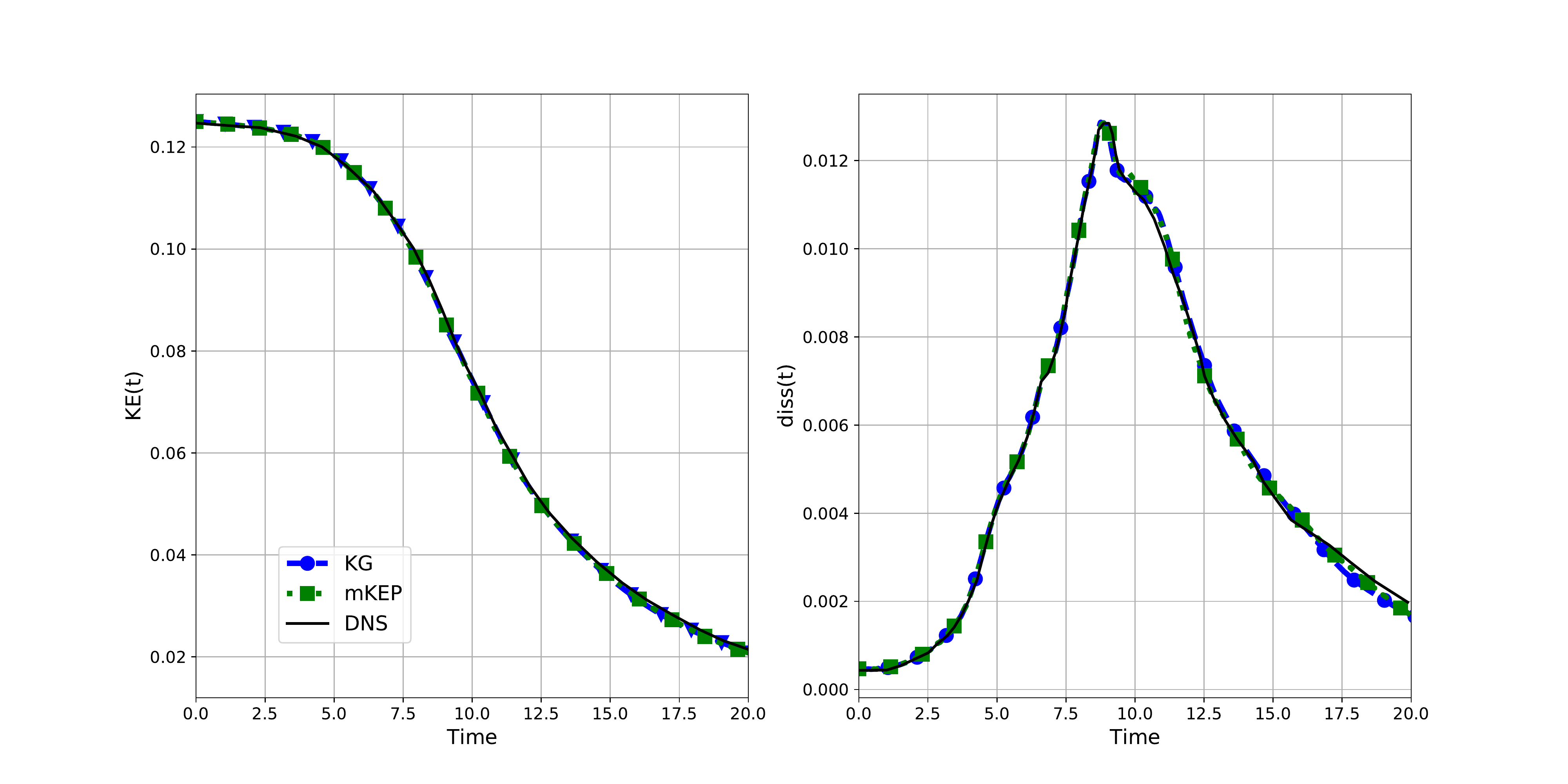}
  \caption{Comparison of kinetic energy and dissipation against DNS of viscous TGV at Re=1600.}
  \label{f:tgv_iles}
  \end{center}
\end{figure}

We use the same initial data as given in Section~\ref{sec:tgv1} and solve the Navier-Stokes equations at a Reynolds number of 1600.  We use a grid consisting of $64^3$ uniform hexahedral elements with $N=3$. We compare the kinetic energy evolution and dissipation rate as a function of time with DNS data from~\cite{DeBonis2013} which is shown in Figure~\ref{f:tgv_iles}; we observe that the new split form performs as well as the KG flux and both of them agree well with the DNS data. 
\subsubsection{Channel flow}
\begin{figure}
  \begin{center}
    \includegraphics[width=0.90\textwidth]{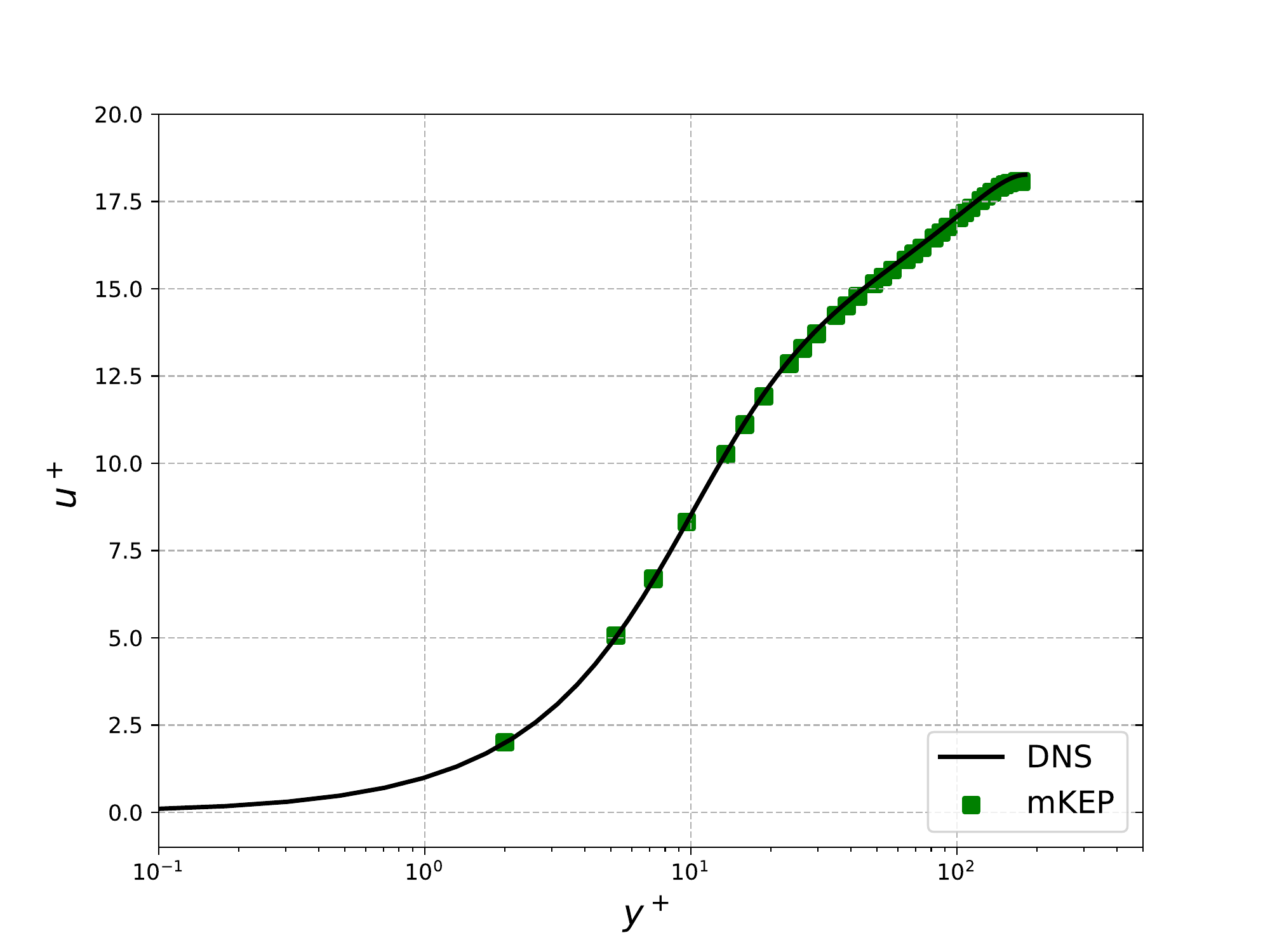}
  \caption{Comparison of mean statistics obtained using mKEP flux with DNS of channel flow for $Re_{\tau} = 180$.}
  \label{f:chan1}
  \end{center}
\end{figure}
\begin{figure}
  \begin{center}
    \includegraphics[width=0.99\textwidth]{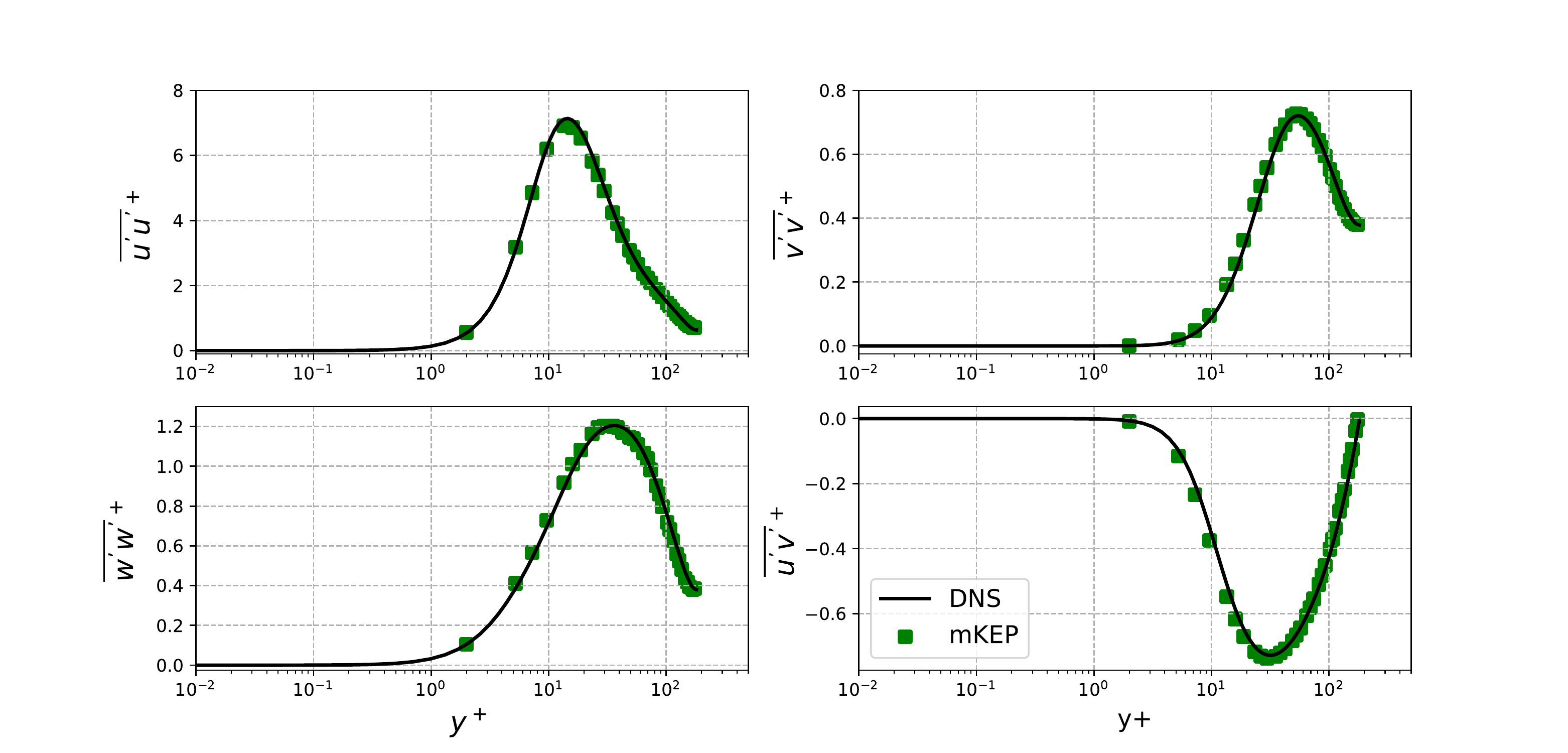}
  \caption{Comparison of variance statistics obtained using mKEP flux with DNS of channel flow for $Re_{\tau} = 180$.}
  \label{f:chan2}
  \end{center}
\end{figure}

To show ILES capabilities for wall bounded flows,  we solve the turbulent channel flow problem~\cite{Kim1987}.  This flow has a periodic domain in the streamwise and spanwise directions with walls at the top and bottom where no-slip condition is imposed.  
The channel has height H and the half height $L=\frac{H}{2}$ decides
the Reynolds' number. 
The flow is driven by a forcing term to get a constant mass flow rate. The test case is characterized by the parameter $Re_{\tau}$ given by
$$
Re_{\tau} = \frac{u_\tau L}{\nu}, \qquad u_\tau = \sqrt{\frac{\tau_w}{\rho}}.
$$

We use a domain length $2\pi$ in the streamwise direction and $\pi$ in the spanwise direction;  the channel height is $H=2$.  It is discretized into $32 \times 22 \times 22$ mesh with parabolic stretching to refine the mesh near the walls.  The computations were compared to DNS data from~\cite{moser_2015} for $Re_\tau=180$. The initial density is taken as $\rho = 1$ and the pressure is chosen such that the Mach number for the initial horizontal velocity $u_0$
is $M = 0.2$.  To this uniform profile, perturbations are added to trigger transition. Figures~\ref{f:chan1} compares the mean streamwise velocity with DNS data, meanwhile Figure~\ref{f:chan2} shows the Reynolds stresses. These figures show excellent agreement of the simulations with DNS data for both first order and second order moments,  showing the suitability of the new split form for ILES simulations. 

\section{Summary and conclusions}

Kinetic energy and entropy preserving split forms are attractive since they can lead to  improved robustness for the computation of turbulent flows, but many such schemes have been shown to be linearly unstable in the sense of eigenvalues when applied to a problem involving the advection of a density wave.  The linearized Euler equations are not strongly stable for this problem since the perturbations can possibly grow with time.  However, a reduced energy is conserved which we use to distinguish the various split form schemes.   We present a new split form for the compressible Euler equations, which we call the modified kinetic energy preserving (mKEP) split form, which conserves the reduced energy.  The central scheme which is stable in terms of eigenvalues does not conserve the reduced energy. For non-normal operators which arise by the linearization of the schemes, eigenvalues characterize only the long time behaviour, but can exhibit energy growth at initial times. We show that our new scheme is consistent with the reduced energy equation and behaves similarly to the Kennedy and Gruber flux in most test cases, while being stable for the density wave problem. Note that the \keep split form also conserves the reduced energy since it has the same linearized equations as the mKEP scheme. To test these schemes, we used a discontinuous Galerkin method which has SBP property.

For the density wave problem, schemes which cannot maintain the constancy of velocity and pressure are always unstable. Schemes which can maintain the constancy like central, Ducros, mKEP and \keep can also be unstable at coarse resolutions due to loss of positivity of density, even though they can be stable in the sense of eigenvalues. At sufficiently high resolutions, such schemes are found to be stable for long times. When the initial condition is perturbed, mKEP and \keep schemes behave similarly and are more stable than central and Ducros schemes. If the initial perturbation is sufficiently large, all schemes fail at large times due to loss of positivity of density, since none of them is monotone.

For the 2D isentropic vortex problem, the \keep and Ducros schemes show better kinetic energy conservation, with Ducros showing a small growth in kinetic energy, while KG and mKEP are sightly more dissipative. The \keep and Ducros are similarly entropy dissipative while the KG and mKEP produce small amounts of entropy. For all schemes, these errors decrease with increasing resolution.

For the 3D invsicid Taylor-Green problem, the \keep scheme shows very high kinetic energy dissipation after the transition to turbulence. The mKEP and KG schemes behave similarly and have small levels of kinetic energy dissipation, while the Ducros scheme has higher dissipation but much less than \keep scheme. The \keep scheme is entropy dissipative while KG, mKEP and Ducros schemes produce small amounts of entropy.

These observations suggest that the mKEP scheme is a good alternative to existing split form schemes when used in a DG method, since it is more stable for the density wave problem, satisfies the reduced energy conservation and has good kinetic energy properties similar to the KG flux. It also performs well for the two viscous turbulent test cases shown here in an ILES setting. The \keep scheme has been shown to perform well in finite difference schemes but it seems to have high levels of kinetic energy dissipation when used in DG methods for turbulent flows. These issues have to be further investigated by applying the schemes to more realistic problems.

\section*{Acknowledgements}
We thank the authors of the Julia code Trixi.jl~\cite{trixi} for making it publicly available. Praveen Chandrashekar was supported by Department of Atomic Energy, Government of India, under  project no. 12-R\&D-TFR-5.01-0520.
\appendix
\section{Symmetrization of linearized equations for density wave problem}
\label{sec:sym}
We attempt to derive an energy equation for the linearized Euler equations~\eqref{eq:lineul} by following the symmetrization approach of~\cite{ABARBANEL1981}. Let $V=[\rho,u,p]$ denote small perturbations and define the change of variables
\[
W = S^{-1} V, \qquad S^{-1} = \begin{bmatrix}
\frac{c}{\sqrt{\gamma} r} & 0 & 0 \\
0 & 1 & 0 \\
-\frac{c}{r\sqrt{\gamma(\gamma-1)}} & 0 & \sqrt{\frac{\gamma}{\gamma-1}} \frac{1}{r c}
\end{bmatrix}, \qquad c = c(x,t) = \sqrt{\frac{\gamma P}{r(x,t)}}
\]
Note that the matrix $S$ is a function of $x,t$ and is not constant. The linearized equations can be written as
\[
W_t + \tilde{A} W_x + \tilde B W = 0, \qquad \tilde{A} = \begin{bmatrix}
V & \frac{c}{\sqrt{\gamma}} & 0 \\
\frac{c}{\sqrt{\gamma}} & V & \sqrt{\frac{\gamma-1}{\gamma}} c \\
0 & \sqrt{\frac{\gamma-1}{\gamma}} c & V
\end{bmatrix}, \qquad
\tilde{B} = \begin{bmatrix}
0 & \frac{c}{r\sqrt{\gamma}} & 0 \\
\frac{c}{2r\sqrt{\gamma}} & 0 & \sqrt{\frac{\gamma-1}{\gamma}} \frac{c}{2r} \\
0 & -\frac{c}{r\sqrt{\gamma(\gamma-1)}} & 0
\end{bmatrix} \pd{r}{x}
\]
If the base flow density is constant, $r =$ constant and $\tilde B=0$, then we obtain a symmetric hyperbolic system with constant coefficients for which the quadratic energy $W^\top W$ is conserved. In the general case, $\tilde A, \tilde B$ are functions of $(x,t)$ and we get the energy equation
\[
W^\top W_t + W^\top \tilde A W_x + W^\top \tilde B W = 0
\]
\[
\half (W^\top W)_t + (W^\top \tilde A W)_x + W^\top \tilde C W = 0
\]
where
\[
\tilde C = \tilde B - \half \tilde A_x = \begin{bmatrix}
0 & \frac{5c}{4r\sqrt{\gamma}} & 0 \\
\frac{3c}{4r\sqrt{\gamma}} & 0 & \sqrt{\frac{\gamma-1}{\gamma}} \frac{3c}{4r} \\
0 & \frac{(\gamma-5)c}{4r\sqrt{\gamma(\gamma-1)}} & 0
\end{bmatrix} \pd{r}{x}
\]
The solution is bounded by
\[
\od{}{t}|W|^2 \le 2 c |W|^2 \qquad\implies\qquad |W(t)| \le \exp(c T) |W(0)|, \qquad 0 \le t \le T
\]
where $c = \norm{\tilde C}$. We are not able to strictly bound the perturbation energy in terms of the initial energy, and the equations do allow for some growth in the solution.
\section{Finite difference schemes}
\label{sec:fd}
In this section, we perform linear stability analysis of some split form fluxes in the finite difference case for the density wave problem.  Consider a finite difference semi-discrete scheme for 1-D Euler equations
\[
\od{U_e}{t} + \frac{F_\eph - F_\emh}{\Delta x} = 0
\]
where $F_\eph = F(Q_e, Q_{e+1})$ is a numerical flux function. We write the numerical flux in terms of the primitive variables $Q=[\rho,u,p]$ which simplifies our analysis. Some examples of numerical fluxes are given in Section~\ref{sec:dg}.
\subsection{Density wave problem in 1-D}
Consider the density wave problem in 1-D introduced in Section~\ref{sec:dwave}, whose exact solution is of the form $\rho(x,t)=r(x,t)=r(x-Vt,0), u=V, p=P$. In this case the Euler flux derivative is of the form $F_x = [ 
\rho_x V, \  \rho_x V^2,  \half \rho_x V^3]^\top$ and all three equations reduce to an advection equation for density. We investigate this property for some numerical schemes.
\paragraph{KG flux}
At time $t=0$, the density and momentum equations reduce to an advection equation
\[
\od{\rho_e}{t}|_{t=0} + V \frac{\avg{\rho}_\eph - \avg{\rho}_\emh}{\Delta x} = 0
\]
which is just the central difference scheme. The energy flux is 
\[
F_E = \frac{1}{\gamma-1} \avg{\rho} \avg{p/\rho} + \avg{p}\avg{u} + \half \avg{\rho} \avg{u^2} \avg{u}
\]
and the energy equation yields
\[
\od{p_e}{t}|_{t=0} + P V \frac{\avg{\rho}_\eph \avg{1/\rho}_\eph - \avg{\rho}_\emh \avg{1/\rho}_\emh}{\Delta x} = 0
\]
The second term on the left is not zero; a Taylor expansion around $x = x_e$ gives
\[
\od{p_e}{t}|_{t=0}+ \frac{PV}{2 \rho_e^3}[ \rho_e \rho'_e \rho''_e - (\rho'_e)^3] \Delta x^3 + O(\Delta x^4) = 0
\]
where primes denote spatial derivatives at $t=0$. While the error in pressure is $O(\Delta x^2)$, its magnitude can be large on coarse grids, expecially in regions where the density $\rho_e$ is small and the pressure and/or velocity are large. When we evolve the solution in time, the scheme changes the pressure distribution in space and thus creates an artificial pressure gradient out of the density gradient, which then changes the velocity also.
\paragraph{Other fluxes.}
A similar effect arises with the KEP flux of Jameson~\cite{Jameson2008} due to the use of enthalpy average, Kuya et al. flux~\cite{Kuya2018}, and all other fluxes in which the pressure and density are coupled in some average. These schemes do not keep the pressure constant which then changes both the velocity and pressure profiles. When we compute the solution using such fluxes, the computations break down due to loss of positivity of solution. The Ducros flux and the modified KEP flux maintain constancy of pressure and velocity.
\subsection{Linearized schemes for density wave problem}
Let $Q_e^0(t) = [ r_e(t), \  V, \  P]^\top$ be solution of the semi-discrete finite volume scheme with constant velocity and pressure. For a scheme to admit such a solution, we require that the momentum and energy fluxes satisfy the following conditions
\begin{equation}
F_m(Q_e^0,Q_{e+1}^0) = P + V F_\rho(Q_e^0,Q_{e+1}^0), \qquad F_E(Q_e^0,Q_{e+1}^0) = \frac{\gamma PV}{\gamma-1} + \half V^2 F_\rho(Q_e^0,Q_{e+1}^0)
\label{eq:fluxpep}
\end{equation}
in which case, all three equations reduce to the density equation. The central, Ducros and modified KEP fluxes satisfy this condition.

Consider a perturbed solution of the form  $Q_e^0 + q_e$ where $Q_e^0$ is the density wave solution and $q_e$ is a small perturbation around this. Let us denote the arguments of the numerical flux by $(X,Y) \to F(X,Y)$ and let $F_X,F_Y$ denote jacobians with respect to the two arguments. Then
\[
F(Q_e^0+q_e, Q_{e+1}^0+q_{e+1}) = F(Q_e^0,Q_{e+1}^0) + F_X(Q_e^0,Q_{e+1}^0) q_e + F_Y(Q_e^0,Q_{e+1}^0) q_{e+1} + O(|q|^2)
\]
The linearized scheme is
\[
\begin{aligned}
\Delta x \begin{bmatrix}
\od{\rho_e}{t} \\
V\od{\rho_e}{t} + \od{}{t}(r_e u_e) \\
\frac{V^2}{2}\od{\rho_e}{t} + V \od{}{t}(r_e u_e) + \frac{1}{\gamma-1}\od{p_e}{t}
\end{bmatrix} + & F_X(Q_e^0,Q_{e+1}^0)q_e + F_Y(Q_e^0,Q_{e+1}^0)q_{e+1} \\
 - & F_X(Q_{e-1}^0,Q_e^0) q_{e-1} - F_Y(Q_{e-1}^0,Q_e^0) q_e = 0
\end{aligned}
\]
Let us also note the following equations for later use
\[
\od{}{t}\left(\half r_e u_e^2 \right) = u_e \od{}{t}(r_e u_e) - \frac{u_e^2}{2} \od{r_e}{t}, \qquad \Delta x \od{r_e}{t} = - V\Delta\avg{r}_e
\]
where $\Delta(\cdot)_e = (\cdot)_\eph - (\cdot)_\emh$.
\subsubsection{Central flux}
The flux jacobian is
{\footnotesize
\[
F_X(Q_e^0,Q_{e+1}^0) = \begin{bmatrix}
\half V & \half r_e & 0 \\ \\
\half V^2 & r_e V & \half \\ \\
\frac{1}{4}V^3 & \frac{\gamma P}{2(\gamma-1)} + \frac{3}{4} r_e V^2 & \frac{\gamma V}{2(\gamma-1)}
\end{bmatrix}, \quad
F_Y(Q_e^0,Q_{e+1}^0) = \begin{bmatrix}
\half V & \half \avg{r} & 0 \\ \\
\half V^2 & \half r_{e+1} V + \half\avg{r}V & \half \\ \\
\frac{1}{4}V^3 & \frac{\gamma P}{2(\gamma-1)} + \half r_{e+1} V^2 + \frac{1}{4}\avg{r}V^2 & \frac{\gamma V}{2(\gamma-1)}
\end{bmatrix}
\]
}
where $\avg{r} = \avg{r}_\eph = \avg{r}_{e,e+1}$. The linearized scheme is given by
\begin{eqnarray*}
\Delta x \od{\rho_e}{t} + V \Delta\avg{\rho}_e + \Delta\avg{r u}_e &=& 0 \\
\Delta x \od{}{t}(r_e u_e) + V \Delta\avg{r u}_e + \Delta\avg{p}_e &=& 0 \\
\Delta x \od{p_e}{t} + V \Delta\avg{p}_e + \gamma P \Delta\avg{u}_e &=& 0
\end{eqnarray*}
From the above equations, we can derive the following relations
\begin{eqnarray*}
\Delta x \od{}{t}\left( \half r_e u_e^2 \right) &=& -V u_e \Delta \avg{ru}_e - u_e \Delta\avg{p}_e + \half u_e^2 V \Delta\avg{r}_e \\
\Delta x \od{}{t}\left(\frac{1}{2\gamma P} p_e^2 \right) &=& -\frac{V}{\gamma P} p_e \Delta\avg{p}_e - p_e \Delta\avg{u}_e
\end{eqnarray*}
We add the above two equations and sum over all cells. Note that
\begin{eqnarray*}
\sum_e p_e \Delta\avg{p}_e &=& \half \sum_e (p_e^2 - p_{e+1}^2) = 0 \\
\sum_e \left[ u_e \Delta\avg{p}_e + p_e \Delta\avg{u}_e \right] &=& \sum_e (p_e u_e - p_{e+1} u_{e+1}) = 0
\end{eqnarray*}
We obtain the energy equation
\[
\Delta x \sum_e \od{\pen_e}{t} = -V\sum_e\left[ u_e \Delta\avg{ru}_e - \half u_e^2 \Delta\avg{r}_e \right] = \frac{V}{2}\sum_e (r_{e+1}-r_e)(u_{e+1} - u_e)^2 \ne 0
\]
We cannot conclude from this that the perturbation energy will be bounded with time.
\subsubsection{Ducros flux}
The flux jacobian is
{\scriptsize
\[
F_X(Q_e^0,Q_{e+1}^0) = \begin{bmatrix}
\half V & \half \avg{r} & 0 \\ \\
\half V^2 & \half r_e V + \half\avg{r}V & \half \\ \\
\frac{1}{4}V^3 & \frac{\gamma P}{2(\gamma-1)} + \half r_e V^2 + \frac{1}{4}\avg{r}V^2 & \frac{\gamma V}{2(\gamma-1)}
\end{bmatrix}, \quad
F_Y(Q_e^0,Q_{e+1}^0) = \begin{bmatrix}
\half V & \half \avg{r} & 0 \\ \\
\half V^2 & \half r_{e+1} V + \half\avg{r}V & \half \\ \\
\frac{1}{4}V^3 & \frac{\gamma P}{2(\gamma-1)} + \half r_{e+1} V^2 + \frac{1}{4}\avg{r}V^2 & \frac{\gamma V}{2(\gamma-1)}
\end{bmatrix}
\]
}
The linearized scheme is given by
\begin{eqnarray*}
\Delta x \od{\rho_e}{t} + V \Delta\avg{\rho}_e + \Delta(\avg{r}\avg{u})_e &=& 0 \\
\Delta x \od{}{t}(r_e u_e) + V \Delta\avg{r u}_e + \Delta\avg{p}_e &=& 0 \\
\Delta x \od{p_e}{t} + V \Delta\avg{p}_e + \gamma P \Delta\avg{u}_e &=& 0
\end{eqnarray*}
The momentum and pressure equations are same as in the case of central flux and hence we obtain the same equation for perturbation energy as in the case of central flux, and this scheme also does not conserve the perturbation energy $\pen$.
\subsubsection{KG flux without advection}
The KG flux does not satisfy the condition~\eqref{eq:fluxpep} required to preserve constant velocity and pressure. It does trivially satisfy this condition if there is no advection, i.e., $V=0$. We will analyze the linear stability in this special case. The flux jacobian is given by
\[
F_X(Q_e^0,Q_{e+1}^0) =  F_Y(Q_e^0,Q_{e+1}^0) = \begin{bmatrix}
0 & \half\avg{r} & 0 \\
0 & 0 & \half \\
0 & \half\alpha P & 0
\end{bmatrix}, \qquad \alpha_\eph = 1 + \frac{\avg{r}_\eph \avg{1/r}_\eph }{\gamma-1}
\]
The linearized equations are given by
\begin{eqnarray*}
\Delta x \od{\rho_e}{t} + \Delta(\avg{r}\avg{u})_e &=& 0 \\
\Delta x \ r_e \od{u_e}{t} + \Delta\avg{p}_e &=& 0 \\
\Delta x \od{p_e}{t} + (\gamma-1) P \Delta(\alpha\avg{u})_e &=& 0
\end{eqnarray*}
Let
\[
\avg{r}\avg{1/r} = 1 + \beta, \qquad\textrm{then}\qquad \alpha = \frac{\gamma}{\gamma-1} + \frac{\beta}{\gamma-1}
\]
and the pressure equation can be written as
\[
\Delta x \od{p_e}{t} + \gamma P \Delta\avg{u}_e + P \Delta(\beta\avg{u})_e = 0
\]
The first two terms are consistent with the linearized pressure equation but the last term is an error term which goes to zero as $\Delta \to 0$ since $\beta \to 0$. However, on any finite mesh, the total energy evolves as
\[
\Delta x \sum_e \od{\pen_e}{t} = - \frac{1}{\gamma} \sum_e p_e \Delta(\beta\avg{u})_e = \frac{1}{\gamma} \sum_e \beta_\eph\avg{u}_\eph(p_{e+1} - p_e) \ne 0
\]
Even in this case of no advection, the KG flux is not consistent with the perturbation equation and  we cannot conclude stability of this scheme since the sign of the right hand side is indeterminate.
\subsubsection{Modified KEP flux}
The flux jacobian is
\[
F_X(Q_e^0,Q_{e+1}^0) = F_Y(Q_e^0,Q_{e+1}^0) = \begin{bmatrix}
\half V & \half \avg{r} & 0 \\ \\
\half V^2 & \avg{r}V & \half \\ \\
\frac{1}{4}V^3 & \frac{\gamma P}{2(\gamma-1)} + \frac{3}{4}\avg{r}V^2 & \frac{\gamma V}{2(\gamma-1)}
\end{bmatrix}
\]
and the linearized scheme is given by
\begin{eqnarray*}
\Delta x \od{\rho_e}{t} + V \Delta\avg{\rho}_e + \Delta(\avg{r}\avg{u})_e &=& 0 \\
\Delta x \od{}{t}(r_e u_e) + V \Delta(\avg{r}\avg{u})_e + \Delta\avg{p}_e &=& 0 \\
\Delta x \od{p_e}{t} + V \Delta\avg{p}_e + \gamma P \Delta\avg{u}_e &=& 0
\end{eqnarray*}
These equations lead to
\begin{eqnarray*}
\Delta x \od{}{t}\left( \half r_e u_e^2 \right) &=& -V u_e \Delta(\avg{r}\avg{u})_e - u_e \Delta\avg{p}_e + \half u_e^2 V \Delta\avg{r}_e \\
\Delta x \od{}{t}\left(\frac{1}{2\gamma P} p_e^2 \right) &=& -\frac{V}{\gamma P} p_e \Delta\avg{p}_e - p_e \Delta\avg{u}_e
\end{eqnarray*}
We add the two equations and sum over all cells to obtain the energy equation
\[
\Delta x \sum_e \od{\pen_e}{t} = -V \sum_e \left[ u_e \Delta(\avg{r}\avg{u})_e - \half u_e^2 \Delta\avg{r}_e \right] = 0
\]
Thus the modified KEP scheme preserves the perturbation energy and is stable in this sense. The density perturbations satisfy the equation
\[
\Delta x \sum_e \od{}{t} \half \rho_e^2 = -V \sum_e \rho_e \avg{\rho}_e - \sum_e \rho_e \Delta(\avg{r}\avg{u})_e = - \sum_e \rho_e \Delta(\avg{r}\avg{u})_e
\]
which is consistent with the continuous equation~\eqref{eq:rhoperttot}, i.e., the convective terms do not lead to any accumulation of density perturbations.
\bibliographystyle{unsrt}  
\bibliography{references}  

\end{document}